\documentclass[12pt]{article}

\usepackage{amstext    }
\usepackage{amsthm    }
\usepackage{a4}
\usepackage[mathscr]{eucal}
\usepackage{mathrsfs}

\usepackage{amsmath}
\usepackage{amssymb}
\usepackage{amscd}

\newtheorem{theorem}{Theorem}[section]
\newtheorem{definition}[theorem]{Definition}
\newtheorem{proposition}[theorem]{Proposition}

\newtheorem{lemma}[theorem]{Lemma}

\newcommand{\cali}[1]{\mathscr{#1}}

\newcommand{\dist}{\mathop{\mathrm{dist}}\nolimits}
\newcommand{\mes}{{\rm mes}}
\newcommand{\Jac}{{\rm Jac}}

\newcommand{\PSH}{{\rm PSH}}

\newcommand{\id}{\mathop{\mathrm{id}}\nolimits}

\newcommand{\Ac}{\cali{A}}

\newcommand{\Cc}{\cali{C}}

\newcommand{\Ec}{\cali{E}}
\newcommand{\Fc}{\cali{F}}

\newcommand{\Hc}{\cali{H}}

\newcommand{\Uc}{\cali{U}}

\newcommand{\C}{\mathbb{C}}

\newcommand{\T}{\mathbb{T}}
\newcommand{\R}{\mathbb{R}}

\renewcommand{\S}{\mathbb{S}}

\renewcommand{\Re}{{\rm Re}}
\renewcommand{\Im}{{\rm Im}}

\begin{document}

\title 
{Corrigendum: Conical plurisubharmonic measure and  new  cross theorems}

\author{Vi{\^e}t-Anh Nguy{\^e}n}

 \maketitle

\begin{abstract} In  the paper  \cite{nv3}  we prove    a theorem on  the boundary behavior of the conical plurisubharmonic  measure.
However,   the proof turns out to be  incomplete.  In the  present work
we give  a corrected  proof of this  theorem.  We  next apply  it to the theory of separately  holomorphic functions.
These applications  are  presented  in a more  accessible  way  than  in   \cite{nv3}.
  \end{abstract}

\noindent
{\bf Classification AMS 2010:}  Primary    32U, 32U05, 32D10, 32D15



\section{Introduction and  statement of the main results}
\label{section_introduction}

Let  $D$ be  an open subset of $\C^n$ and  $A\subset\partial D.$
We suppose in addition that
 $D$ is  {\it locally $\Cc^2$ smooth} on $A$ (i.e.
for any $\zeta\in A,$ there exist an open neighborhood $U=U_\zeta$ of
 $\zeta$ in $\C^n$ and a   real function
 $\rho=\rho_\zeta\in \Cc^2(U)$ such that $D\cap U=\lbrace z\in U:\
 \rho(z)<0\rbrace$ and $d\rho(\zeta)\not=0$).
For $\zeta\in  A$ and  $1<\alpha<\infty,$ we consider the conical  approach  region
 \begin{equation*}
 \Ac_\alpha(\zeta):=\left\lbrace  z\in D:\  \vert  z-\zeta\vert <\alpha\cdot \dist(z,T_\zeta)  \right\rbrace,
 \end{equation*}
where  $ \dist(z,T_\zeta)$ denotes the Euclidean    distance  from the  point $z$ to the tangent hyperplane $T_\zeta$  of $\partial D$ at $\zeta.$

For any function $u$ defined on $D,$ let
\begin{equation*}
 \hat{u}(z):=
\begin{cases}
u(z),
  & z\in   D,\\
\sup\limits_{\alpha>1} \limsup\limits_{w\in \Ac_\alpha(z),\ w\to z}u(w), & z \in \partial D.
\end{cases}
\end{equation*}
 Next, consider the function $h_{A,D}:=\sup\limits_{u\in\mathcal{F}}u,$ where
\begin{equation*}
\Fc:=\left\lbrace   u\in\PSH(D):\ u\leq 1\ \text{on}\ D,\
   \hat{u}\leq 0\ \text{on}\ A    \right\rbrace.
\end{equation*}
Here $\PSH(D)$ denotes the set   of all functions  plurisubharmonic
on $D.$
Then
 the {\it conical plurisubharmonic measure} of $A$ relative to $D$ is given
 by
\begin{equation}
 \omega(z,A,D):=
  h^{\ast}_{A,D}(z)    , \qquad z \in  D,
\end{equation}
where $u^{\ast}$ denotes the  upper semicontinuous regularization of a function $u.$

 A manifold  $M\subset\C^n$ of class $\Cc^2$   is  said to be  {\it generic} if, for every point
$z\in M,$ the complex linear  hull of the tangent space  $T_zM $ (to $M$  at  $z$)
coincides with the  whole space $\C^n.$

The  main purpose  of this  work  is  to investigate  the boundary behavior of the conical plurisubharmonic  measure in  a   special  but  important case, and  thereafter to apply this  study to the  theory of separately holomorphic mappings.
     Now  we  are in the position to state  the main result.


\begin{theorem}\label{main_thm}
  Let $M\subset \C^n$ be  a generic  manifold of  class  $\Cc^2$  and  $D$  a domain
  in $\C^n$  such that  $M\subset \partial D$ and  $D$ is locally $\Cc^2$ smooth on $M.$ Let $A\subset M$  be a  measurable subset  of positive  measure.
  Then for all density points $z$ relative to $A,$
$\widehat{\omega}(z,A,D)=0.$
\end{theorem}

This   theorem  describes the  stable character of the  the conical plurisubharmonic  measure   $\omega(\cdot,A,D) $ along the conical  approach regions at   all density points
relative to  $A.$ It sharpens  the previous  results  of A. Sadullaev (see \cite{sa}) and B. Coupet (see \cite{co})
where  the estimate  $\omega(\cdot,A,D)<1$ on $D$ was  obtained.
Our proof relies  on  the  use of families of analytic  discs attached to $M$   and  on some fine
estimates of plurisubharmonic  functions.

\smallskip

This paper is organized as follows.

\smallskip

We begin Section 2 by
collecting some results of the  method  of attaching analytic  discs
  to a  generic  manifold in the spirit of Coupet's work  \cite{co}.  Next,  we develop necessary estimates for the conical plurisubharmonic measure
and then prove  Theorem  \ref{main_thm}.
Section 3  concludes  the  article  with  various  applications
of  Theorem  \ref{main_thm} in the  theory of separately holomorphic mappings.

\medskip

\indent{\it{\bf Acknowledgment and  comments.}}
The first version of  the paper    has been published  in \cite{nv3}, but fortunately in March 2011         Malgorzata Zajecka  (Krak\'ow)  found   a  gap therein. Namely,   the author claimed    in the proof of Theorem 2.1  in  \cite{nv3}
that  this  theorem should follow implicitly  from Th\'eor\`eme 2 in \cite{co}.  However, this claim is  not correct.
 The   present version  
 has  filled   this gap.   But the main  idea of the proof  is  always  the same as in  \cite{nv3}.
More specifically,   we   first  construct   analytic  discs attached to a given generic manifold,  and  then apply     
estimates for plurisubharmonic  functions.    In the   construction  of  analytic  discs Theorem 2.1 in   \cite{nv3}  has  to be  replaced by  Proposition \ref{prop_Coupet}  and \ref{prop_VA} below. Since  the   subsequent steps     rely on this  construction, 
the  details  of the proof  are different from  those   given in \cite{nv3}.
 However,
  the  strategy  as well  as the       lemmas    are almost  unchanged in comparison with    \cite{nv3}.

The  author would like to thank  Malgorzata Zajecka  for the  valuable help.

\section{Proof}
For  $x\in\R^m$  let $|x|$  denotes  its Euclidean norm. 
For $x\in\R^n$ and $r>0$ let $B(x,r)$ be the Euclidean ball with center $x$ and with radius $r.$ 
For  a   $\Cc^2$ smooth Riemannian manifold $M$  of dimension $m$ let $\mes_M$ denote the
$m$-dimensional Lebesgue measure on $M$. When there is no fear of confusion we often write
$\mes$ instead of   $\mes_M.$  If $M$ is  a smooth submanifold  in $\R^m$ then  we   equip  $M$  with the  Riemannian metric
induced from $\R^m.$  
For  two    functions  $A$ and $B,$ we use  the  following   conventional notation.   We  write $A\lesssim B$  (or equivalently 
$A=O(B)$)
  if  there is  a constant $c>0$ such that  $|A|\leq c|B|.$
We write  $A\approx B$ if $A\lesssim B$ and $B\lesssim A.$
Moreover, $A=o(B)$ means that  $ { |A|\over |B|}\to 0$ as $|B|\to 0.$
 For  a differentiable map $g:\ M\to N$  between Riemannian manifolds,
 let $\Jac g$ denote the Jacobian  matrix of $g.$ If, moreover, $\dim M=\dim N$ then   we denote  by $|\Jac g|$  the determinant of $\Jac g.$ 

A smooth  generic  manifold $M\subset\C^n$  is  said to be  {\it totally real} if  $\dim_\R M=n.$
Our proof  will be  divided  into two  cases.
In the first one  we assume that $M$ is   totally real.   The second one  will  treat the general case  of  $M.$  

To deal  with the  first  case,  let $M\subset \C^n$ be  a  totally manifold of class $\Cc^2.$ We  may   assume without loss of generality that  $0\in M$ and
$T_0M=\R^n$  (it suffices to perform  an affine change of coordinates). 
In the  sequel  given  $z\in \C^n$ we often  write     $z=x+iy$ with  $x,y\in\R^n.$  $M$  is  then  defined   in a neighborhood
of $0\in \C^n$  by the  equation $z=x+ih(x),$ where $h$ is  a function of class $\Cc^2$  defined  in an open  neighborhood  of $0\in \R^n$   with  values in $\R^n$    satisfying  $h(0)=0$ and  $dh(0)=0.$

  A holomorphic  disc on a Jordan domain $\Sigma$  is, by definition, a  continuous   map $f:\  \overline{\Sigma}\rightarrow\C^n$
such that $f|_\Sigma$ is  holomorphic.    A holomorphic  disc $f$ is  said  to be {\it attached to}  $M$ on  an arc $\Gamma\subset \partial\Sigma$ if
  $f(\Gamma)\subset M.$
  
 Let $\Delta$  be the open unit disc in $\C$  and  $\T:=\partial \Delta.$ Let $T$ be the  conjugate  operator on $L^2(\T)$, that is, the operator which associates  to  every $u\in L^2(\T)$ an element
$T(u)\in   L^2(\T)$  such that $\int_\T T(u)=0$ and that $u+iT(u)$ is  the boundary value of a holomorphic function on $\Delta.$
Let $\tilde{T}(u)$ be the harmonic  extension of $T(u)$ on $\Delta.$

Fix a smooth function $\phi$   defined on $\overline{\Delta}$ harmonic on $\Delta$
such that
\begin{itemize}
\item[$\bullet$] $\phi=0$ on $\{e^{i\theta}\in\T:\ \vert \theta\vert \leq \frac{\pi}{2}\}; $
\item[$\bullet$]  
$\phi<0$ on   $ \{e^{i\theta}\in\T:\    {\pi\over 2}<\vert \theta\vert \leq \pi\}    ;$ 
\item[$\bullet$]  $\max_{\theta\in[-\pi,\pi]} |T(\phi)(e^{i\theta})|=1;$
 \item[$\bullet$] $\frac{\partial \phi(se^{i\theta})}{\partial s}|_{s=1,\theta=0}>0.$
 \end{itemize}
Here  $[0,1]\times [-\pi,\pi]\ni (s,\theta)\mapsto se^{i\theta}\in\overline{\Delta}.$

For  $u=(u_1,\ldots,u_n)\in L^2(\T)^n $ let  $T(u)$ denote
$(T(u_1),\ldots,T(u_n)).$  

Consider  the following Pinchuk's version of  Bishop equation
$$x=-T(h\circ x)+\zeta-T(\phi)\tau.$$  
 As  was  shown in  \cite{co}, there  exists $\delta>0$ such that the  above  equation admits a unique  solution  $x\in L^2(\T)^n$  when $(\zeta,\tau)\in\R^{2n}$  with $|\zeta|<\delta$
and  $|\tau|<\delta$  are  given.   
By the uniqueness of  the solution   of the  above equation, 
we get that
\begin{equation}\label{eq_special_soultion}
 x(e^{i\theta},\zeta,0)=\zeta,\qquad  \theta\in\R,\ \zeta\in\R^n, \  |\zeta|<\delta.
 \end{equation}
Consider  the following  map $g:\  \T\times \{ \zeta\in\R^n:\   |\zeta|<\delta\}\times  \{ \tau\in\R^n:\   |\tau|<\delta\}\rightarrow\R^{2n}$ defined by
\begin{equation}\label{eq_map_g}
g(w,\zeta,\tau):= x(w,\zeta,\tau)+i\Big(h\circ x (w,\zeta,\tau)+\tau  \phi (w)\tau\Big).
\end{equation}
Let  $\tilde{g}$ be  the map defined  on $ \Delta\times \{ \zeta\in\R^n:\   |\zeta|<\delta\}\times  \{ \tau\in\R^n:\   |\tau|<\delta\}$
with values in $\R^{2n}$  such that
      $\tilde{g}(\cdot,\zeta,\tau)$  is  holomorphic  on $\Delta$ and admits  $g(\cdot,\zeta,\tau)$ as its boundary value.
 In particular, we deduce  from  (\ref{eq_special_soultion}) and (\ref{eq_map_g}) that 
\begin{equation}\label{eq_g(w,0,0)}
  \tilde{g}(w,0,0)=0,\qquad  w\in\Delta.
  \end{equation}

The  following result  will be  needed.

\begin{proposition}\label{prop_Coupet}   Let $1<\alpha<\infty.$
\begin{itemize}
\item[(i)] Let $ \Theta:\  \overline{\Delta}\times \big(\Ac_\alpha(0)\cap B(0,\delta)\big)\rightarrow\C^n$ be  given by
$$\Theta(w,\zeta,\tau)=\tilde{g}(w,\zeta,\tau)-\zeta+T(\phi)(w) \cdot\tau -i\phi(w)\cdot\tau.$$
Then
  $|  \Theta(w,\zeta,\tau)| \lesssim  |\tau|^2,$   ${\partial \Theta(w,\zeta,\tau)\over \partial \zeta}\lesssim |\tau|,$
 ${\partial \Theta(w,\zeta,\tau)\over \partial \tau}\lesssim |\tau|,$
${\partial \Theta(se^{i\theta},\zeta,\tau)\over \partial \theta}\lesssim |\tau|^2.$

\item[(ii)]   There  exist
 a  subdomain $\Sigma\subset \Delta$ of the form
$$
\Sigma:=\{ z=se^{i\theta},\   s_0\leq s <1,\ \theta\in[-{\pi\over 2},{\pi\over 2}]:  e^{i\theta}\in \Gamma\},
$$
with  $\Gamma\subset \T$  an open arc  containing $1$ and $0<s_0<1$ a constant,  
and  an open neighborhoods $U$  of $0$ in $B(0,\delta)$
such that 
$$
\tilde{g}(w,\zeta,\tau)\in  D,\qquad (w,\zeta,\tau)\in \overline{\Sigma}\times (\Ac_\alpha(0)\cap U).
$$
\end{itemize}
\end{proposition}
\begin{proof}
It follows  implicitly from  Th\'eor\`eme 2 in  Coupet's work \cite{co}. For the sake of clarity we  recall briefly  his  argument.

 First recall from   estimates  (1) and (2) in \cite{co}  that
\begin {eqnarray*}
{\partial \tilde{g}\over\partial \zeta}(w,0,0)&=&\id_\zeta,\\
{\partial \tilde{g}\over\partial \tau}(w,0,0)&=&-T(\phi)(w)\id_\tau+ i\phi(w)\id_\tau.
\end{eqnarray*}
This, combined  with (\ref{eq_g(w,0,0)}) and the formula for $\Theta,$  implies that   
$$\Theta(w,\zeta,\tau)=\tilde{g}(w,\zeta,\tau)-\tilde{g}(w,0,0)-{\partial \tilde{g}\over\partial \zeta}(w,0,0)(\zeta)-{\partial \tilde{g}\over\partial \tau}(w,0,0)(\tau). $$
On the other hand,  for  $(\zeta,\tau)\in \Ac_\alpha(0)\cap B(0,\delta),$  we  have, by the definition of conical approach regions,
  that
$$
|(\zeta,\tau')|<\alpha \tau_n,
$$
where $\tau=(\tau',\tau_n).$  So 
\begin{equation*}
|\zeta|<\alpha|\tau|.
\end{equation*}
This, together with the last identity for $\Theta,$ imply  all the  estimates  of Part (i) except the last one.
   Arguing  as  in the proof of   estimate  (3) in   \cite{co} and using the last  estimate  we can show 
that, for $  (\zeta,\tau)\in  \Ac_\alpha(0) \cap B(0,\delta),$
\begin{equation}\label{eq_radial_derivatives}
\begin{split}
{\partial \Im \tilde{g}\over \partial s}(se^{i\theta},\zeta,\tau)&={\partial  \tilde{\phi}\over \partial s}(se^{i\theta})\cdot \tau +o(\tau),
\\
{\partial \Re \tilde{g}\over \partial s}(se^{i\theta},\zeta,\tau)&={-\partial  \tilde{T}(\phi)\over \partial s}(se^{i\theta})\cdot \tau +o(\tau).
\end{split}
\end{equation}
This  implies that
\begin{equation*}
{\partial  \tilde{g}\over \partial s}(se^{i\theta},\zeta,\tau)={\partial\over\partial s}\Big[\zeta- T(\phi)(w) \cdot\tau+i \phi(w)\tau\Big]+O(|\tau|^2).
\end{equation*}
Hence, 
$$
{\partial \Theta(w,\zeta,\tau)\over \partial s}=O(|\tau|^2).
$$
On the other hand, since   $\Theta(\cdot,\zeta,\tau)$ is holomorphic in $w,$  it follows that
 $$
{\partial \Theta(w,\zeta,\tau)\over \partial \theta}=   is {\partial \Theta(w,\zeta,\tau)\over \partial s}   .
$$
This, coupled  with the last  inequality, gives the last estimate of   Part (i).

Fix   $\alpha',$  $\alpha''$ such that  $1<\alpha<\alpha'<\alpha''<\infty.$  Since $D$ is  locally  $\Cc^2$-smooth on $A\ni 0,$
we may find   an open neighborhood $U$ of $0$ in $\C^n$ such that
\begin{equation}\label{eq_conic_approach}
\begin{split}
&\quad \Ac_{\alpha''}(z)\cap U\subset D,\qquad  z\in \partial D\cap U,\\
&\quad  y+z \in  \Ac_{\alpha''}(z),\qquad  z\in \partial D\cap U,\ y\in \Ac_{\alpha'}(0)\cap U.
\end{split}
\end{equation}
 On the  other hand,  
estimate (\ref{eq_radial_derivatives}) for $s=1$ gives  that  
\begin{eqnarray*}
{\partial \Im \tilde{g}\over \partial s}(e^{i\theta},\zeta,\tau)&=&{\partial  \tilde{\phi}\over \partial s}(e^{i\theta})\cdot \tau +o(\tau),\\
{\partial \Re \tilde{g}\over \partial s}(e^{i\theta},\zeta,\tau)&=&   o(\tau),
\end{eqnarray*}
where  the last  estimate  follows  from  the  identity
$${\partial  \tilde{T}(\phi)\over \partial s}(e^{i\theta})=-{\partial  \phi\over \partial \theta}(e^{i\theta})=0,\qquad |\theta|<{\pi\over 2}.
$$
Consequently, by shrinking $U$ if necessary we can choose  a  subdomain $\Sigma\subset \Delta$ of the form stated  by  Part (ii)
 such that the ray emanating from  $\tilde{g}(e^{i\theta},\zeta,\tau)$ and passing through  $\tilde{g}(se^{i\theta},\zeta,\tau)$
 cuts  the unit sphere $\S$ (in $\R^{2n}$)  at a point $\eta$ such that  $\dist(\eta,\nu)=O(|\tau|),$
 where 
$$\nu:={ \big(0,{\partial  \tilde{\phi}\over \partial s}(e^{i\theta})\cdot \tau\big) \over  \big|{\partial  \tilde{\phi}\over \partial s}(e^{i\theta})\cdot \tau
 \big|}\in\S.
 $$
  Moreover,  by shrinking $\Sigma$ if necessary, we  may assume that ${\partial  \tilde{\phi}\over \partial s}(e^{i\theta})>0$ for $e^{i\theta}\in \Gamma.$ 
 Since this  assumption implies that $\nu\in \Ac_\alpha(0),$ it follows that $\tilde{g}(se^{i\theta},\zeta,\tau)-\tilde{g}(e^{i\theta},\zeta,\tau)\in  \Ac_{\alpha'}(0)$
 for $  (\zeta,\tau)\in  \Ac_\alpha(0) \cap B(0,\delta)$ and $\delta>0$  small  enough.
So, by (\ref{eq_conic_approach}),
$$ 
 \tilde{g}(se^{i\theta},\zeta,\tau)\in  \Ac_{\alpha'}\big( \tilde{g}(e^{i\theta},\zeta,\tau)\big)\cap U\subset D,\qquad  se^{i\theta}\in\Sigma.
$$
 This completes Part (ii).
\end{proof}

 Following  Th\'eor\`eme 2 in Coupet's work \cite{co}
 we  will    construct a map $G$  using $g$ and $\tilde{g}.$ This  new map  is the  key ingredient  for our proof of  Theorem \ref{main_thm}. 
 Consider the map  $(\zeta,\tau)\mapsto g(1,\zeta,\tau)$ from  an open neighborhood of $0$ in $\R^{2n}$
 to $M.$  We know from  estimates (1) and (2) in \cite{co} that the rank of   this map is $n$ at $0\in\R^{2n}.$
Since $g(1,0,0)=0,$  it follows that there is a  $\Cc^1$ map $a$ defined  on an open neighborhood of $0$ in $\R^n$
with values in $\R^n$ such that
\begin{equation*}
a(0)=0,\quad \text{and}\quad  g(1,\zeta,\tau)=0\ \text{if and only if}\ \zeta=a(\tau).
\end{equation*}         
 Fix a   $\Cc^\infty$ map $b$  defined on $\R^n\setminus \{0\}$ to the space of linear  endomorphisms of rank $n-1$ from $\R^{n-1}$ to $\R^n$   such that  
 $b(\tau)(\R^{n-1})$ is orthogonal to $\tau$ and that $b(\tau)=b(r\tau)$, $r>0.$ Consequently, by shrinking $\Sigma$ if necessary, we may find  an open neighborhood $U$ of $0\in\R^{2n}$
such that  
  the following map is  well-defined
 \begin{equation}\label{eq_G}
G(w,\zeta',\tau):= \tilde{g}(w,a(\tau)+b(\tau)(\zeta'),\tau),\qquad  w\in \Sigma,\ (\zeta',\tau)\in \Hc , 
\end{equation}
where
$$\Hc=\Hc(U,\alpha):=  \left\{ (\zeta',\tau)\in \R^{n-1}\times \R^n:\  |\zeta'|\leq 4 |\tau|\  \&\  (0,\tau)\in\Ac_\alpha(0)\cap U \right \}$$
In what follows,   for $r>0$ and $z\in M$ let $B_M(z,r):=B(z,r)\cap M.$   Some  important properties of  $G$  are  listed  below.
  
 \begin{proposition}\label{prop_VA}
 Under the  above hypothesis and notation
 \begin{itemize}
\item[(i)]  For $(\zeta',\tau)\in\Hc,$ $G(\cdot,\zeta',\tau)$ is  a holomorphic disc on $\Sigma$
 with image    contained  in 
 $D.$
Moreover,   $G(\cdot,\zeta',\tau)$ is attached to $M$
on $\Gamma$   and  $G(\cdot,0,0)=0$ on $\overline{\Sigma}.$
\item[(ii)] 
For each fixed $s_0\leq s\leq  1,$ consider  the map $G_s:\  \Gamma\times\Hc\rightarrow\C^n$ defined by
$G_s(w,\zeta',\tau):= G(sw,\zeta',\tau).$  Then, for $s_0\leq s<  1,$  
 $G_s$ maps  $ \Gamma\times\Hc$ diffeomorphically  onto  its  image.
 Moreover,     $|\Jac G_s(w,\zeta',\tau)|\approx (1-s)^n|\tau|$  as $\tau\to 0.$ 
\item[(iii)]  There is  a constant $c>1$ such that
for every fixed $\tau\in\R^n$ with $(0,\tau)\in \Ac_\alpha(0)\cap U$ and every ${1\over 4}\leq\beta\leq 4,$ $G_1(\cdot,\cdot,\tau)$ maps $\Gamma\times\{\zeta'\in \R^{n-1}:\  |\zeta'|\leq  \beta|\tau|\}$ diffeomorphically
onto  an open neighborhood $V$ of $G_1(1,0,\tau)$ in $M$  such that
 $$  B_M\big(G_1(1,0,\tau),{1\over c}|\tau|\big)\subset V\subset  B_M\big(G_1(1,0,\tau),c|\tau|\big).$$
 Moreover,
  $| \Jac G_1(\cdot,\cdot,\tau)|\approx |\tau|.$    
\end{itemize}
\end{proposition}
\begin{proof} 
Part (i) follows  from identity (\ref{eq_G}) and Part (ii) of  Proposition \ref{prop_Coupet}. 
Next,
writing $G_s=\Re G_s +i \Im G_s$ with
$  \Re G_s=((\Re G_s)_1,\ldots, (\Re G_s)_n),$  $  \Im G_s=((\Im G_s)_1,\ldots, (\Im G_s)_n),$
and  $w=se^{i\theta},$
we have  that 
\begin{equation*}
\Jac G_s =\left(
\begin{array}{ccccccc}
{\partial (\Re G_s)_1\over\partial \theta} &  {\partial (\Re G_s)_1\over\partial \zeta_1} &\cdots  & {\partial (\Re G_s)_1\over\partial \zeta_{n-1}}
&{\partial (\Re G_s)_1\over\partial \tau_1} &\cdots  & {\partial (\Re G_s)_1\over\partial \tau_n} \\
\vdots&  \vdots  & \vdots &\vdots &\vdots &\vdots  &\vdots \\
{\partial (\Re G_s)_n\over\partial \theta} &  {\partial (\Re G_s)_n\over\partial \zeta_1} &\cdots  & {\partial (\Re G_s)_n\over\partial \zeta_{n-1}}
&{\partial (\Re G_s)_n\over\partial \tau_1} &\cdots  &{ \partial (\Re G_s)_n\over\partial \tau_n} \\
{\partial (\Im G_s)_1\over\partial \theta} &  {\partial (\Im G_s)_1\over\partial \zeta_1} &\cdots  & {\partial (\Im G_s)_1\over\partial \zeta_{n-1}}
&{\partial (\Im G_s)_1\over\partial \tau_1} &\cdots  &{ \partial (\Im G_s)_1\over\partial \tau_n} \\
\vdots&  \vdots  & \vdots &\vdots &\vdots &\vdots  &\vdots \\
{\partial (\Im G_s)_n\over\partial \theta} &  {\partial (\Im G_s)_n\over\partial \zeta_1} &\cdots  & {\partial (\Im G_s)_n\over\partial \zeta_{n-1}}
&{\partial (\Im G_s)_n\over\partial \tau_1} &\cdots  & {\partial (\Im G_s)_n\over\partial \tau_n} \\
\end{array}
\right)
\end{equation*}
To compute $ \Jac G_s,$ 
using  Part (i) of Proposition  \ref{prop_Coupet} and  identity (\ref{eq_G}) we have that
  \begin{equation}\label{eq_G_new}
 G(w,\zeta',\tau)=a(\tau)+b(\tau)(\zeta')-T(\phi)(w)\tau-i\phi(w)\cdot\tau+ \Theta(w, a(\tau)+b(\tau)(\zeta')     ,\tau).
 \end{equation}
 Using this  together with  the  inequalities  obtained  in  Part (i) of Proposition  \ref{prop_Coupet},
 we  can    estimate  the  entries  of $\Jac G_s$ as  well as $\Jac G_1(\cdot,\cdot,\tau).$  
More specifically,
a straightforward computation gives  the following  estimates  on the entries of $\Jac G_s$
\begin{equation}\label{eq_estimates_entries}
\begin{split}
{\partial (\Re G_s)_k\over\partial \theta}&=
-{\partial T(\phi)( w)\over\partial \theta}\cdot \tau_k + O(|\tau|^2),\qquad  1\leq k\leq n;\\
{\partial (\Re G_s)_k\over\partial \zeta_j}&=
 O(1),\qquad  1\leq k\leq n,\ 1\leq j\leq n-1;\\
{\partial (\Re G_s)_k\over\partial \tau_j}&=
O(1),\qquad  1\leq k\leq n,\ 1\leq j\leq n  ;\\
 {\partial (\Im G_s)_k\over\partial \theta}&=
 { \partial \phi(w)\over\partial\theta}\cdot\tau_k+ O(|\tau|^2),  \qquad  1\leq k\leq n           ;\\
{\partial (\Im G_s)_j\over\partial \tau_j}&=
 \phi(w)+ O(|\tau|),  \qquad  1\leq j\leq n;\\
{\partial (\Im G_s)_k\over\partial \tau_j}&=
  O(|\tau|),  \qquad  1\leq j\not=k\leq n;\\
{\partial (\Im G_s)_k\over\partial \zeta_j}&=
  O(|\tau|),  \qquad  1\leq j,k\leq n.
\end{split}
\end{equation}
 Using the properties of  $\phi$  we  see that the quantities $|{\partial T(\phi)( w)\over\partial \theta}|$
and $|{ \partial \phi(w)\over\partial\theta}|$ are $\approx 1$  and that  $-\phi(w)\approx  1-s.$
On the  other hand,
applying Part (i) of Proposition \ref{prop_Coupet} to (\ref{eq_G_new}) yields that
$$
{\partial  \Re G_s(w,\zeta',\tau)\over \partial \zeta'}=  b(\tau)+O(|\tau|).
$$
This, combined with  the fact that  $b(\tau)$ is  orthogonal to $\tau$ and 
the first  estimate  in (\ref{eq_estimates_entries}), implies that   
$$ |\Jac \Re G_s(\cdot,\cdot,\tau)|\approx |{\partial \Re G_s\over\partial \theta}|\approx |\tau|. $$
Now  we prove Part (iii).
When $s=1,$    we know from Part (i) that  $G_1(\cdot,\cdot,\tau)$  maps $\Gamma\times  \{\zeta'\in\R^{n-1}:\  |\zeta'|\leq 2|\tau|\}$
to $M.$ Observe  that  $T_0M$ is spanned  by $\Re G_s.$ So $\Jac G_1(\cdot,\cdot,\tau)= \Jac \Re G_s(\cdot,\cdot,\tau) +O(|\tau|).$
Using the  above  discussion on $\Jac \Re G_s(\cdot,\cdot,\tau),$   the inclusions of Part (iii) follows.
Moreover, 
  $$|\Jac G_1(\cdot,\cdot,\tau)|= |\Jac \Re G_s(\cdot,\cdot,\tau)|\approx |\tau|,$$
 which   proves the last assertion  of Part (iii).

To prove Part (ii) observe from   (\ref{eq_estimates_entries})  that
$$
|\Jac \Im G_s(w,\zeta',\cdot)|\approx \phi(w)^n=  (1-s)^n.
$$ 
On the other hand, expanding  the determinant of $ \Jac G_s$ by fixing its  first  $n$  rows
and  using     (\ref{eq_estimates_entries}) again,  we  see that  the dominant  element in 
the resulting expression   is  $$ |\Jac \Re G_s(\cdot,\cdot,\tau)|\cdot |\Jac \Im G_s(w,\zeta',\cdot)|\approx |\tau|(1-s)^n.$$
This  finishes the proof of Part (ii). 
 \end{proof}

 We need the following elementary lemma.
 \begin{lemma}\label{lem_ele}
 Let $\Gamma$ be  an open arc in $\T\cap \partial\Sigma$ of the form  $\Gamma=\{ e^{i\theta}:\ \theta\in (-2\theta_0,2\theta_0)\},$  where $0<\theta_0<\frac{\pi}{2}$ is a fixed number.    Then  there exists
  $1<c<\infty$
 with the following  property. If   $u$ is a subharmonic  function  defined in $\Sigma$ with $u\leq 1$ on $\Sigma$
 and $\hat{u}\leq 0$ on $B,$ where  $B$ is  a measurable subset of $ \Gamma$ with $\frac{\mes(\Gamma\setminus B)}{\mes(\Gamma)}\leq \epsilon^2,$
  then $\sup_{\rho=1-\epsilon,\ |\zeta| < \theta_0} u(\rho e^{i\zeta})\leq  c\epsilon.$
 \end{lemma}
\begin{proof}
  Let $\psi:\ \Sigma\to \Delta$ be   a conformal map with $\psi(1)=1.$
  Since $\psi$ is  smooth  on $\partial\Sigma\cap\overline{\Delta},$  we may reduce  the problem  to  the  case 
where $\Sigma=\Delta.$
 The assumption on $u$ implies  that for $\rho\in [0,1]$ and $\vartheta\in [-\pi,\pi],$
\begin{equation*}
u(\rho e^{i\vartheta})\leq  \frac{1}{2\pi} \int\limits_{-\pi}^{\pi} \frac{(1-\rho^2) 1_{B}}{\vert  e^{i\theta} -\rho e^{i\vartheta}\vert^2}d\theta,
\end{equation*}
where
\begin{equation*}
 1_B(\theta):=
\begin{cases}
0,
  & e^{i\theta} \in   B,\\
 1, & e^{i\theta}\not \in B.
\end{cases}
\end{equation*}
Therefore, an easy estimate  shows that for $\rho=1-\epsilon$ and $\vert\vartheta\vert<\theta_0.$
\begin{equation*}
u(\rho e^{i\vartheta})\leq \frac{1}{2\pi}\int\limits_{-2\theta_0\epsilon^2}^{2\theta_0\epsilon^2}
\frac{(1-\rho^2)}{\vert  1 -\rho\vert^2}d\theta +\frac{1}{2\pi} \int\limits_{\vert \theta\vert >2\theta_0} \frac{(1-\rho^2) }{\vert  e^{i\theta} -\rho e^{i\vartheta}\vert^2}d\theta\leq   c\epsilon.
\end{equation*}
\end{proof}

Suppose without loss of generality that $0$ is  a  point of density relative to $A$ in $M.$
The proof of Theorem  \ref{main_thm}  is  divided into  several steps.
In the first  three steps   we assume that $M$ is  totally real.  Fix   arbitrary  $0<\epsilon_0<1$ and $1<\alpha_0<\infty.$    We  only need  to show that there exists
a  sufficiently small open neighborhood  $\Uc$ of $0\in \C^n$   such that
    $u(z)\leq \epsilon_0$  for all
   $z\in \Ac_{\alpha_0}(0)\cap \Uc$  and   $u\in \Fc,$  where  the family $\Fc$  was   defined  in Section  \ref{section_introduction}.

\noindent {\bf Step 1:}  {\it  Construction of a subset $\Omega_r\subset D\cap B(0,r)$ for all $r>0$ small enough such that
$u(z)\leq  \epsilon_0,$  $z\in\Omega_r,$  $u\in \Fc.$ }

The idea  is to use  families of holomorphic discs  attached to $M$  which parametrize a   open neighborhood of $0$ in $M.$ These families are supplied by Proposition \ref{prop_VA}.

  For $r>0$ let $U_r:=\{ \tau\in\R^n:\  {r\over 2}\leq |\tau|\leq r\}$ and $V_r:=\{ \zeta'\in\R^{n-1}:\   |\zeta'|\leq 2 r\}.$
 Let $\gamma>1$ and $\alpha >\alpha_0$ be  constants large enough whose  exact value will be  determined  later on.
Since $0$ is  a point of density relative to  $A,$  we may find $0<r_0\ll 1$ such that  for all $0<r<r_0,$
\begin{equation}\label{eq0_main_thm}
\frac{\mes\big(B_M(0,r)\setminus A\big)}{ \mes\big(  B_M(0,r)\big)}<\epsilon_0^\gamma.
\end{equation}

Fix  $0<r<r_0$ and  $\tau\in U_r.$
By Part (iii) of Proposition \ref{prop_VA},   
$G(\cdot,\cdot,\tau)$ maps $\Gamma\times V_r$ diffeomorphically
onto  an open neighborhood  $M_\tau$ of $0$ in $M$  which is  $\approx B_M(0,r) .$
 Moreover,
  $ \Jac G(\cdot,\cdot,\tau)\approx r.$   
This, combined  with  (\ref{eq0_main_thm}), implies that
\begin{equation*} 
\frac{\mes\big(  M_\tau\setminus A     \big)}{\mes(M_\tau)}<\epsilon_0^\gamma.
\end{equation*}
Since  $M_\tau$ is parametrized by the  family of real  curves $\{ G(\Gamma,\zeta',\tau):\  \zeta'\in V_r\},$  it follows from the last estimate that
\begin{equation}\label{eq1_main_thm}
\frac{\mes\big( V_r\setminus  V^\tau_r  \big)}{\mes\big( V_r \big)     }<\epsilon_0^{\gamma\over 2},
\end{equation}
where
$$
 V_r^\tau:=\left\lbrace \zeta'\in V_r:\ {  \mes\big(  G(\Gamma,\zeta',\tau)\setminus A\big)  \over
   \mes\big( G(\Gamma,\zeta',\tau)\big) }< \epsilon_0^{\gamma\over 2}   \right\rbrace,\qquad  \tau\in U_r.
 $$
 For $0<\epsilon<1$  let
  $\Gamma_\epsilon:= \{(1-\epsilon)e^{i\theta}:\ \theta\in (-\theta_0,\theta_0)\}.$
   Applying Lemma \ref{lem_ele} to the holomorphic discs $G(\cdot,\zeta,\tau)|_\Sigma$  attached to $M$
on $\Gamma$ which are supplied by  Proposition  \ref{prop_VA}, we deduce from  (\ref{eq1_main_thm}) that
 \begin{equation}\label{eq1bis_main_thm}
 u(z)\leq  c\epsilon_0,\qquad  z\in \Omega_r,\ u\in\Fc,
 \end{equation}
 where 
 $$\Omega_r:=   \bigcup\limits_{  \tau\in U_r,\  \zeta'\in V_r^\tau   }G(\Gamma_{\epsilon_0},\zeta',\tau)   .$$
 Note that  by  Part (i) of Proposition \ref{prop_Coupet} and  the  above formula  for $\Omega_r,$ we get that 
 $ \Omega_r\subset D\cap B(0,r).$
 
\noindent {\bf Step 2:}  {\it  For $0<r<r_0,$ let $ H_r:=   \left\lbrace z=x+iy\in\C^n:\ 
 z\in \Ac_\alpha(0),\  y_n\in (r,2r)  \right\rbrace.$  Then 
$$ {\mes (H_r\setminus \Omega_r)\over\mes(H_r)} \leq \epsilon_0^{{\gamma\over 2}-n+1}.$$}

Consider the set
 $$
 R_r:=\bigcup\limits_{ \tau\in U_r,\  \zeta'\in V_r} G(\Gamma_{\epsilon_0},\zeta',\tau)
 $$
 Since  $V_r^\tau\subset V_r$ we  clearly have that $\Omega_r\subset R_r.$  Moreover, we deduce  from (\ref{eq1_main_thm}) that
 \begin{equation}\label{eq2_main_thm}
{\mes (R_r\setminus \Omega_r)\over  \mes (R_r)}<\epsilon_0^{\gamma\over 2}.
 \end{equation}
  Next  we will introduce  a set $S_r$ such that that  $\Ac_\alpha(0)\cap R_r\subset S_r \subset R_r.$
  For $  \theta\in (-\theta_0,\theta_0)$ and $   \tau\in U_r,$
  let
  $
 V_{\theta,\tau}$ be the  set of all $ \zeta'\in V_r$ such that
$$ |\Re G((1-\epsilon_0)e^{i\theta},\zeta',\tau)|\leq  \alpha \epsilon_0r.
$$
To  estimate  the quotient ${\mes(V_{\theta,\tau})\over  \mes (V_r)}$ 
observe that
\begin{eqnarray*}
&& |\Re G((1-\epsilon_0)e^{i\theta},\zeta',\tau) - \Re G(1,\zeta',\tau)|\\
&=&
 |\Re\tilde{g} ((1-\epsilon_0)e^{i\theta}, a(\tau)+b(\tau)(\zeta'),\tau)
 -\Re\tilde{g} (1, a(\tau)+b(\tau)(\zeta'),\tau)\\
 &\leq &\big(T(\phi)((1-\epsilon_0)e^{i\theta})-T(\phi)(1)\big)\tau+O(|\tau|^2)\\
&\lesssim & \epsilon_0r,
\end{eqnarray*}
where the first inequality follows  from Part (i) of Proposition \cite{co}.
Moreover, since  $\tilde{g}(1,a(\tau) ,\tau)=0,$ we  get that
 \begin{eqnarray*}
 |\Re G(1,\zeta',\tau)|&=&|\Re \tilde{g}(1,a(\tau)+b(\tau)(\zeta') ,\tau)|\\
 &=&|\Re \tilde{g}(1,a(\tau)+b(\tau)(\zeta') ,\tau)-\Re \tilde{g}(1,a(\tau) ,\tau)|
 =b(\tau)(\zeta')+O(|\tau|^2),
\end{eqnarray*}
where the last equality also follows  from Part (i) of Proposition \cite{co}.
This  together  with the  last observation  imply that
$$
 |\Re G((1-\epsilon_0)e^{i\theta},\zeta',\tau) -b(\tau)(\zeta')|  \lesssim  \epsilon_0r+O(|\tau|^2).
$$
 Hence,  we conclude,  for $\alpha\gg 1,$ that $
 V_{\theta,\tau}$ is approximatively  a ball  with center $0$ and radius $\epsilon_0r$ in $V_r.$
Consequently, 
$${\mes(V_{\theta,\tau})\over  \mes (V_r)}  \approx \epsilon_0^{n-1}.$$
 Define
$$
 S_r:=\left  \{G(  \Gamma_0,\zeta',\tau):\        \tau\in U_r,\  \zeta'\in V_{\theta,\tau}\right\}.
 $$
 Since  $V_{\theta,\tau}\subset V_r,$  we  clearly have that $S_r\subset R_r.$ Moreover, using Part (ii) of Proposition \ref{prop_VA}  the last estimate 
 gives that
  \begin{equation}\label{eq3_main_thm}
 {\mes(S_r)\over \mes (R_r)}\approx {\mes(V_{\theta,\tau})\over  \mes (V_r)}  \approx \epsilon_0^{n-1}.
 \end{equation}
 Using Part (i) of Proposition \ref{prop_Coupet}   and the estimate 
   \begin{equation*}
   - \phi( (1-\epsilon_0)e^{i\theta})\approx \epsilon_0,\qquad  \theta\in (-\theta_0,\theta_0),
   \end{equation*}
 we  see easily that
  $   \Ac_\alpha(0)\cap B(0,r) \subset S_r.$  Hence, 
$  \Ac_\alpha(0)\cap B(0,r) \subset S_r\subset R_r.$ 
  This,  combined with  (\ref{eq2_main_thm}), implies  that
  $$
  {\mes \big ( ( \Ac_\alpha(0)\cap B(0,r)      \setminus \Omega_r)\over  \mes (R_r)}\leq
{\mes  ( S_r      \setminus \Omega_r)\over  \mes (R_r)}\leq
   {\mes  ( R_r      \setminus \Omega_r)\over  \mes (R_r)} \leq \epsilon_0^{\gamma\over 2}.
$$
 Using this and    (\ref{eq3_main_thm})  and the estimate  $ {\mes( \Ac_\alpha(0)\cap B(0,r))\over  \mes( S_r)}\approx 1$
 it follows  that
$$
{\mes \big ( ( \Ac_\alpha(0)\cap B(0,r))      \setminus \Omega_r)\over  \mes (  \Ac_\alpha(0)\cap B(0,r))}
\lesssim  \epsilon_0^{{\gamma\over 2}-n+1}.
$$
Since $H_r\subset   \Ac_\alpha(0)\cap B(0,r)$ and   $ {\mes( \Ac_\alpha(0)\cap B(0,r))\over  \mes( H_r)}\approx 1$
the desired conclusion of Step 2 follows.

\noindent{\bf  Step 3:} {\it  End  of the proof of  Theorem  \ref{main_thm}  when $M$ is  totally real.}

 First  we  need  some elementary lemmas.

\begin{lemma}\label{lem1}
For   $0<a,b<\infty,$    there  exists a    constant $ c>0$ that depends only on the quotient ${a\over b} $
  with the  following property.
  Consider the  domains
\begin{eqnarray*}
I&:=&\left\lbrace z=x+iy\in\C^n:\ x_1,\ldots,x_n,y_1,\ldots,y_{n-1}\in (-2b,2b),\ 
   y_n\in (\frac{a}{2},4a)  \right\rbrace,\\
J&:=&\left\lbrace z=x+iy\in\C^n:\ x_1,\ldots,x_n,y_1,\ldots,y_{n-1} \in  (-b,b),\  y_n\in (a,2a)  \right\rbrace.
\end{eqnarray*}
Then, for  every  $u\in \PSH(I)$ and every $0<\epsilon<1$ such that $u\leq  1$ on $I$
and that
\begin{equation*}
 \frac{\mes\big( \{z\in I:\  u(z)\geq  \frac{\epsilon}{2}\} \big)}
 { \mes( I  )      }  < c\epsilon,
\end{equation*}
we have  $u<\epsilon$ on $J.$
\end{lemma}
\begin{proof} Observe that there  exists the maximum number $0<r<\infty$ dependent only on $a$ and $b$  such that
 the ball $B(z,r)$ centered at $z$ with radius $r$ in $\C^n$ is contained in $I$ for all $z\in \overline{J}.$
 By the sub-mean property of plurisubharmonic functions we have
 \begin{equation*}
 u(z)\leq {1\over \mes(B(z,r))} \int_{B(z,r)}u(w)dw,\qquad z\in J.
 \end{equation*}
 Setting $c:={\mes(B(z,r))\over 2\mes(I)},$ we see that $c$ depends only on  ${a\over b}.$ 
 Moreover,  we have that for every $z\in J,$
 $$
 { \mes(\big( \{z\in B(z,r):\  u(z)\geq  {\epsilon\over 2}\} \big)\over  \mes(B(z,r))}\leq
{\mes(I)\over  \mes(B(z,r))}\cdot { \mes(\big( \{z\in I:\  u(z)\geq  {\epsilon\over 2}\} \big)\over  \mes(I)}     
 \leq {\epsilon\over2}.
 $$
  This, combined with the above sub-mean estimate,  implies    that
$$
u(z)\leq  \big(1-{\epsilon\over 2}\big){\epsilon\over 2}+{\epsilon\over 2}<\epsilon, \qquad  z\in J.
$$
\end{proof}

\begin{lemma}\label{lem2}
For  every $\alpha>1$  
there  exists a    constant  $ c>0$   with the  following property.
For  all $a>0$    consider the  domains
\begin{eqnarray*}
 H&:=&\left\lbrace z=x+iy\in\C^n:\   z\in \Ac_{ 8\sqrt{n}\alpha}(0),\  y_n\in (\frac{a}{2},4a)  \right\rbrace,\\
K&:=&\left\lbrace z=x+iy\in\C^n:\  z\in \Ac_\alpha(0),\  y_n\in (a,2a)  \right\rbrace.
 \end{eqnarray*}
Then for  every  $u\in \PSH(H)$ and every $0<\epsilon<1$   such that $u\leq  1$ on $H$
and that
\begin{equation*}
 {\mes\big( \{z\in H:\  u(z)\geq  \frac{\epsilon}{2}\} \big)  \over \mes(  H )      }  \leq c\epsilon,
\end{equation*}
we have  $u<\epsilon$ on $K.$
 \end{lemma}
 \begin{proof}
 Applying  Lemma  \ref{lem1}  to the  case  where $b:=2\alpha  a$  and  observing that
 $I\subset  H,$  $K\subset  J,$ the desired conclusion  follows.
 \end{proof}

  To prove Theorem  \ref{main_thm} in the case  where $M$ is  totally real,
 let $0<r<r_0$ be  an arbitrary  number.  Choose  $\alpha:=8\sqrt{n}\alpha_0.$  
By Step 2, we get that
\begin{equation*}
 u(z)\leq  c\epsilon_0,\qquad  z\in \Omega_r,\ u\in\Fc,
 \end{equation*}
 where $\Omega_r$ is  a subset of $D$  that satisfies 
$$ {\mes (H_r\setminus \Omega_r)\over\mes(H_r)} \leq \epsilon_0^{{\gamma\over 2}-n+1}.$$
Now it suffices  to choose  $\gamma\geq 2n+1.$ Then the above quotient is  dominated by $c\epsilon_0$ when $\epsilon_0$ is  small enough.
So we are in the position to apply
  Lemma  \ref{lem2}. Consequently,  
  $u(z)\leq \epsilon_0$ for all $z=x+iy\in \Ac_{\alpha_0}(0)$ with $  y_n\in (r,2r).$
Since $0<r\ll 1$ is arbitrary, this  completes the proof.

\noindent {\bf Step 4:}  {\it The general case.}

Let $M$  be a generic manifold of dimension $n+k$  in $\C^n.$  By a  complex linear change of coordinates
we may  assume that  $T_0M=\R^{n-k}\times \C^k.$ 
Let  $U(k)$ be the  group of  unitary matrices  of  rank $k.$  For  every $H\in U(k),$ there is  a totally real  manifold
$0\ni M^H\subset M$  such that $T_0M^H:=     \R^{n-k}\times H(\R^k).$
  Observe that  $M^H\subset M\subset \partial D.$ 
Since $0$ is  a point of density relative to  $A,$  a slicing  argument  shows that  we may find $0<r_0\ll 1$   
such that  for every $0<r<r_0,$
there exists an $H\in U(k)$  such that 
\begin{equation*}
\frac{\mes\big(B_{M^H}(0,r)\setminus A\big)}{ \mes\big(  B_{M^H}(0,r)\big)}<\epsilon_0^\gamma.
\end{equation*}
Using this  instead of (\ref{eq0_main_thm}),  we  argue  as in the  proof of Steps 1, 2 and  3  above.   Consequently, the conclusion of Step 3 and  the above  observation together 
imply that
$$ \lim_{z\to 0,\  z\in\Ac_\alpha(0)} u(z)=0,\qquad u\in \Fc,\  1<\alpha <\infty.$$
Hence, the proof of the last  step is  finished.
  \hfill $\square$

\section{Applications}  \label{section_applications}


%
%
%
The purpose of this  section is   to  derive from Theorem A in  our previous work \cite{nv1} and  Theorem \ref{main_thm}   a  boundary cross theorem  and a mixed cross theorem. We  first recall  some terminology  and notation (relative to the  system  of conical approach  regions) introduced in Section  2  of    \cite{nv1}.

\begin{definition} Let  $D\subset \C^m$ be  an open set  and   $A\subset\overline{D}.$   
We say that   $A$ is {\rm
locally pluriregular} at a point $a\in \overline{A}$  if one  of the  following   cases happens.
\begin{itemize}
\item[$\bullet$]  $a\in D$ and 
$\omega(a,A\cap U, U)=0$
  for  all open neighborhoods $U$ of $a$ in $D,$ where $\omega(\cdot, A\cap U, U)$  is  the Siciak relative extremal function. 
  \item[$\bullet$]  $a\in \partial D$ and $D$ is  locally $\Cc^2$ smooth  on $\{a\},$ and 
$\widehat{\omega}(a,A\cap\partial D\cap U,D\cap U)=0$
  for  all open neighborhoods $U$ of $a,$  where $\omega(\cdot, A\cap\partial D\cap U,D\cap U)$  is  the  conical plurisubharmonic measure. 
  \end{itemize}
Moreover, $ A$ is said to be {\rm locally pluriregular } if it is locally
pluriregular at all points $a\in A.$ 
\end{definition}
 
 In what follows let $M\subset \C^m$ be  a generic  manifold of  class  $\Cc^2$  and  $D$  an open set
  in $\C^m$  such that  $M\subset \partial D$ and that $D$ is locally $\Cc^2$ smooth on $M.$ 
 Now we are  in the position to formulate the following  version
of the  plurisubharmonic measure.
\begin{definition} \label{defi_conical_measure}
Let $M\subset \C^m$ be  a generic  manifold of  class  $\Cc^2$  and  $D$  a open set
  in $\C^m$  such that  $M\subset \partial D$ and  $D$ is locally $\Cc^2$ smooth on $M.$ Let $A\subset M.$ 
 Let 
$\widetilde{A}  :=\bigcup\limits_{P\in \Ec(A)} P,$ where
\begin{equation*}
\Ec(A):=\left\lbrace P\subset M:\  P\ \text{is  locally pluriregular,}\
\overline{P}\subset  A       \right\rbrace,
\end{equation*}
 The {\rm  generalized conical plurisubharmonic measure of $A$  relative to $D$} is
 the function $\widetilde{\omega}(\cdot,A,D)$
 defined by
\begin{equation*}
\widetilde{\omega}(z,A,D):=  \omega(z,\widetilde{A},D),\qquad  z\in D.
\end{equation*}
\end{definition}
The following result  is  a  consequence of Theorem \ref{main_thm}.
\begin{lemma}\label{lem_consequence}
Let $A$ be a  measurable subset of $M$  such that   $\mes_M(A)>0.$  Then 
\begin{equation*}
\widetilde{\omega}(z,A,D)\leq  \omega(z,A',D),\qquad  z\in D,
\end{equation*}
where  $A'$ is the set of density points   of $A.$
\end{lemma} 
  \begin{proof}
  By Theorem \ref{main_thm},    $A$  is  locally  pluriregular  at all  points
   of $A'.$ On the other hand, since $\mes_M(A\setminus A')=0,$  every point of $A'$ is  also  a density point relative to $A'.$
    Hence,   $A'$  is  locally pluriregular.
 Choose an increasing sequence
$(A_k)_{k=1}^{\infty}$  of subsets of $A$  such that
  $A_k\ $  is  closed and  $
\mes_M \Big( A\setminus \bigcup\limits_{k=1}^{\infty} A_n\Big)=0.$ 
Observe that  $A'_k$ is locally  pluriregular,  $\overline{A'}_k \subset A_k\subset A$ 
       Hence,  by Definition \ref{defi_conical_measure}  we get that $\widehat{A}:= \bigcup\limits_{k=1}^{\infty} A'_k\subset
 \Ec(A).$ 
 So
\begin{equation*}
\widetilde{\omega}(z,A,D)\leq  \omega(z,\widehat{A},D),\qquad  z\in D.
\end{equation*}
 Since  $\mes_M  (A'\setminus \widehat{A})=0$ 
we deduce  from Theorem \ref{main_thm} that
\begin{equation*}
\omega(z,\widehat{A},D)\leq  \omega(z,A',D),\qquad  z\in D.
\end{equation*}
This, combined with the previous  estimate, implies the lemma.
\end{proof}

Now  we are in the position  to give two applications of Theorem \ref{main_thm}.
The first one is  a boundary cross theorem.
  \begin{theorem} \label{application1}
  Let $M\subset \C^m$  and $N\subset \C^n$ be  two generic  manifolds of  class  $\Cc^2.$  Let $D\subset \C^m$
and $G\subset \C^n$ be two open sets  such that  $M\subset \partial D,$ $N\subset\partial G$ 
and  that $D$ is locally $\Cc^2$ smooth on $M,$ $G$ is locally $\Cc^2$ smooth on $N.$
 Let $A\subset M,$ $B\subset N$ be  measurable  sets  such that $\mes_M(A)>0,$  $\mes_N(B)>0.$  
   Define
   \begin{eqnarray*}
   W&:=& \big ((D\cup A)\times B\big)\bigcup  \big(  A\times(G\cup B)    \big) ,\\
    \widehat{W'} & := &\left\lbrace  (z,w)\in D\times G:\  \omega(z,A',D)+\omega(w,B',G)<1     \right\rbrace,
   \end{eqnarray*}
     where $A'$  (resp.  $B'$) is  the set of  all density  points   of $A$  (resp.  $B$).

 Let   $f:\     W   \longrightarrow \C$ be a function
  which satisfies the   following  conditions:
   \begin{itemize}
   \item[$\bullet$]   For every $a\in A,$  $f(a,\cdot)$   is  holomorphic 
on $G$ and  continuous  on $G\cup B.$ For every $b\in B,$  $f(\cdot,b)$   is  holomorphic 
on $D$ and  continuous  on $D\cup A.$
 
    \item[$\bullet$] $f$ is locally bounded, that is,  for every $x\in W,$ there exist  an open neighborhood $U_x$ of $x$ in $W$
and a finite constant $c_x$   such that $\sup_{U_x} |f|< c_x.$
    \item[$\bullet$]          $f|_{A\times B}$is    continuous.
    \end{itemize}
    Then there exists a  unique holomorphic function
$\hat{f}$ defined  on $\widehat{W'}$  such that  for all $  1<\alpha<\infty,$
  \begin{eqnarray*}
\lim\limits_{(z,w)\to (\zeta,\eta),\  z\in \Ac_\alpha(\zeta),\ w\in \Ac_\alpha(\eta)} \hat{f}(z,w)
 &=&f(\zeta,\eta),\qquad
  (\zeta,\eta)\in  A'\times B';\\
\lim\limits_{(z,w)\to (\zeta,\eta),\  z\in \Ac_\alpha(\zeta)} \hat{f}(z,w)
 &=&f(\zeta,\eta),\qquad
  (\zeta,\eta)\in  A'\times G;\\
\lim\limits_{(z,w)\to (\zeta,\eta),\   w\in \Ac_\alpha(\eta)} \hat{f}(z,w)
 &=&f(\zeta,\eta),\qquad
  (\zeta,\eta)\in  D\times B'.
\end{eqnarray*}
 If,  moreover,  $\sup_W\vert f\vert<\infty,$ then
\begin{equation*}
 \vert \hat{f}(z,w)\vert\leq \big(\sup_{A\times B}\vert f\vert\big)^{1-\omega(z,A',D)-\omega(w,B',G)} \big (\sup_W\vert
 f\vert\big)^{\omega(z,A',D)+\omega(w,B',G)},\qquad (z,w)\in\widehat{W'}.
\end{equation*}
\end{theorem}
\begin{proof} 
 Combing Theorem A in \cite{nv1}  and    Lemma \ref{lem_consequence},  the theorem follows.
\end{proof}

The second  application is  a mixed cross   theorem.

 \begin{theorem} \label{application2}
  Let $M\subset \C^m$    be  	a generic  manifolds of  class  $\Cc^2.$  Let $D\subset \C^m$
 be an open set  such that  $M\subset \partial D$ 
and  that $D$ is locally $\Cc^2$ smooth on $M.$ 
 Let $A\subset M$  be a measurable  set  such that $\mes_M(A)>0.$  
 Let $ G\subset \C^n$ be  an open set  and let  $B$  be  a locally pluriregular subset of
  $ G$.
    
   Define
    \begin{eqnarray*}
   W  &:= &    \big ((D\cup A)\times B\big)\bigcup  \big(  A\times(G\cup B)    \big)     ,\\
    \widehat{W'}  &:= &\left\lbrace  (z,w)\in D\times G:\  \omega(z,A',D)+\omega(w,B,G)<1     \right\rbrace,
   \end{eqnarray*}
    where $A'$  is the set of  density points  of $A.$

Let   $f:\ W\longrightarrow \C$  be  a function 
  which satisfies  the  following conditions:
   \begin{itemize}
   \item[$\bullet$]    For every $a\in A,$  $f(a,\cdot)$   is  holomorphic 
on $G.$   For every $b\in B,$  $f(\cdot,b)$   is  holomorphic 
on $D$ and  continuous  on $D\cup A.$
    \item[$\bullet$] $f$ is  locally bounded  along  $A\times G.$
    \end{itemize}
    Then  there  exists a unique holomorphic function 
$\hat{f}$ defined on $\widehat{W'}$ 
such that  
  \begin{eqnarray*}
 \lim\limits_{(z,w)\to (\zeta,\eta),\  z\in \Ac_\alpha(\zeta)} \hat{f}(z,w)
 &=&f(\zeta,\eta),\qquad
  (\zeta,\eta)\in  A'\times G,\  1<\alpha<\infty ;\\
\lim\limits_{(z,w)\to (\zeta,\eta) } \hat{f}(z,w)
 &=&f(\zeta,\eta),\qquad
  (\zeta,\eta)\in  D\times B'.
\end{eqnarray*}
 If,  moreover,  $\sup_W\vert f\vert<\infty,$ then
\begin{equation*}
 \vert \hat{f}(z,w)\vert\leq \big(\sup_{A\times B}\vert f\vert\big)^{1-\omega(z,A',D)-\omega(w,B,G)} \big (\sup_W\vert
 f\vert\big)^{\omega(z,A',D)+\omega(w,B,G)},\qquad (z,w)\in\widehat{W'}.
\end{equation*}
 \end{theorem}
\begin{proof} 
 Combing Theorem A in \cite{nv1}  and    Lemma \ref{lem_consequence},  the theorem follows.
\end{proof}

Observe that  Theorem \ref{application1} (resp.  Theorem \ref{application2})  is  a particular
case    of Theorem 10.4 (resp. Theorem  10.5) in \cite{nv1}. However, the proofs of the latter theorems
are analogous  to those of the former ones.

 Theorem \ref{main_thm} also  implies  Corollary 2 and 3  in \cite{nv2} which  were stated without proof therein.
These corollaries generalize Theorem  \ref{application1} and \ref{application2}  to the case  where
some pluripolar  or thin singularities are allowed  (see \cite{jp1} or \cite{np} for more details
on this  issue). We refer the reader to  the  book  by Jarnicki and Pflug   \cite{jp2} for a comprehensive and  systematic  introduction to the theory of
separately holomorphic functions. 

\noindent
V.-A.  Nguy{\^e}n,
Vietnamese Academy  of Science  and  Technology,
Institute of Mathematics,
Department  of Analysis,
18  Hoang Quoc  Viet  Road, Cau Giay  District,
10307 Hanoi, Vietnam.\\
{\tt nvanh@math.ac.vn}

\noindent
{\sc Current  address:}
Math{\'e}matique-B{\^a}timent 425, UMR 8628, Universit{\'e} Paris-Sud,
91405 Orsay, France.\\
  {\tt VietAnh.Nguyen@math.u-psud.fr}, {\tt http://www.math.u-psud.fr/$\sim$vietanh}

\end{document}